\documentclass[a4paper]{article}
\usepackage[T1]{fontenc}
\usepackage{geometry}
\usepackage{amssymb, amsmath, color}
\usepackage[amsmath, thmmarks]{ntheorem}
\theoremheaderfont{\bfseries\upshape} 
\theoremseparator{ :}
\newtheorem*{de}{Definition}
\newtheorem*{theo}{Theorem}
\newtheorem{cl}{Claim}

\theorembodyfont{\normalfont}
\newtheorem*{rem}{Remark}
\theoremsymbol{\ensuremath{\square}}
\newtheorem*{pr}{Proof}
\title{Maximization of the second conformal eigenvalue of spheres}
\author{Romain Petrides \footnote{Romain Petrides, UMPA-ENS Lyon 46 allée d'Italie 69364 Lyon Cedex 07.
romain.petrides@ens-lyon.fr}}
\date{}
\begin{document}
\maketitle

\begin{abstract} We establish in this paper an upper bound on the second eigenvalue of $n$-dimensional spheres in the conformal class of the round sphere. This upper bound holds in all dimensions and is asymptotically sharp as the dimension increases. 
\end{abstract}

Given $(M,g)$ a smooth compact Riemannian manifold (without boundary), the spectrum of the Laplacian $\Delta_g = -div_g\left(\nabla_g\right)$ is a discrete sequence of eigenvalues 
$$0=\lambda_0\left(M,g\right)< \lambda_1\left(M,g\right)\le \lambda_2\left(M,g\right)\le \dots \le \lambda_k\left(M,g\right)\le \dots$$
which goes to $+\infty$ as $k\to +\infty$. The two first eigenvalues are simple, the eigenfunctions associated to $\lambda_0=0$ being the constant functions. A natural, and often adressed, question is to get estimates on the eigenvalues thanks to some geometric assumptions. In this paper, we discuss maximisation of eigenvalues for metrics in a given conformal class with fixed volume. We focus on the case of the standard sphere. 
 
\medskip We let $\mathbf{S}^n$ be the unit sphere of $\mathbf{R}^{n+1}$ for $n\ge 2$. If $g$ is a metric on $\mathbf{S}^n$, we are interested in the scale invariant quantity 
$$\Lambda_{n,k} (g)=\lambda_k(\mathbf{S}^n,g)Vol_g(\mathbf{S}^n)^{\frac{2}{n}}$$
In dimension $2$, we can maximize $\Lambda_{2,k}$ on regular metrics. An inequality has been proved for $k=1$ by Hersch~\cite{HER1970} :
$$ \Lambda_{2,1}(g) \leq 8\pi $$
with equality iff $g$ is the round metric. He followed the proof of the maximization by Szegö~\cite{SZE1954} of the first non zero Neumann eigenvalue for planar domains, attained by discs. Nadirashvili found an optimal maximization for $k=2$. He proved in \cite{NAD2002} that 
$$ \Lambda_{2,2}(g) < 16\pi $$
where the supremum is attained in the degenerate case of the union of two identical spheres. His idea was  used later in \cite{GIR2009} to show that among simply connected planar domains, the second non zero Neumann eigenvalue is maximal in the degenerate case of two discs of the same area.

\medskip If we look for an analogous inequality in dimension $n\geq 3$, we have to restrict our attention to some classes of metrics since $\Lambda_{n,k}$ is not bounded on the set of regular metrics (see \cite{COL1994}). It is natural, as suggested in~\cite{ELS1986} and~\cite{COL2003}, to consider the set of metrics in some conformal class. Indeed, in any given conformal class, $\Lambda_{n,k} (g)$ admits some upper bound (see \cite{KOR1993}). Thus we define the conformal spectrum of $\left(\mathbf{S}^n,[g_0]\right)$, where $[g_0]$ is the class of metrics conformal to the round metric $g_0$, by
$$ \lambda_k^c(\mathbf{S}^n,[g_0]) = \sup_{g\in[g_0]} \Lambda_{n,k}(g) $$
The theorem of Hersch was generalized in this framework in \cite{ELS1986}. We have that 
$$ \lambda_1^c(\mathbf{S}^n,[g_0]) = n \sigma_n^{\frac{2}{n}} $$
where $\sigma_n$ is the volume of the unit $n$-dimensional sphere. We know almost nothing about $\lambda_k^c(\mathbf{S}^n,[g_0])$ for $k\geq 2$. A lower bound was obtained by a method of conformal surgery in \cite{COL2003}. For all $k$, we have that 
$$ \lambda_k^c(\mathbf{S}^n,[g_0]) \geq n (k \sigma_n)^{\frac{2}{n}}\hskip.1cm. $$
Nadirashvili, Girouard and Polterovich conjectured in \cite{GIR2009} that this inequality is an equality in all dimensions for $k=2$, where the supremum is attained for the union of two identical spheres~:

\medskip\noindent {\bf Conjecture (\cite{GIR2009})}~: {\it for any metric $g\in \left[g_0\right]$, 
$$\lambda_2\left(\mathbf{S}^n,g\right) Vol_g\left(\mathbf{S}^n\right)^{\frac{2}{n}} < n\left(2\sigma_n\right)^{\frac{2}{n}}\hskip.1cm.$$
}

\medskip In the way to this conjecture, the following theorem gives an "asymptotically sharp" upper bound :

\begin{theo} Let $n\geq 2$ and $g\in \left[g_0\right]$ a metric on $S^n$ conformal to the round metric. Then 
$$\lambda_2\left(\mathbf{S}^n , g\right)Vol_g\left(\mathbf{S}^n\right)^{\frac{2}{n}}  < K_n n(2\sigma_n)^{\frac{2}{n}}$$
where $K_n$ is a constant independant of $g\in[g_0]$ given by 
$$ K_n = \frac{n+1}{n} \left( \frac{\Gamma(n)\Gamma(\frac{n+1}{2})}{\Gamma(n+\frac{1}{2})\Gamma(\frac{n}{2})} \right)^{\frac{2}{n}}\hskip.1cm.$$
\end{theo}

Note that $K_2=1$, that $1< K_n \le 1.04$ for all $n\ge 3$ and that ${\displaystyle \lim_{n \to \infty} K_n = 1}$. The theorem is sharp in dimension $2$ and was in fact already proved by Nadirashvili in \cite{NAD2002}. In \cite{GIR2009}, Girouard, Nadirashvili and Polterovich established this inequality in odd dimensions. 

\medskip We prove in this paper this theorem in all dimensions, unifying the previous proofs in dimension $n=2$ and in odd dimensions and by the way extending it. The starting point of the proof is a construction, described in section \ref{construction} below, initiated by Nadirashvili \cite{NAD2002} and used by Girouard, Nadirashvili and Polterovich \cite{GIR2009} in odd dimension. However, our use of this construction differs from that of these two papers : we use the min-max characterisation of the second eigenvalue up to the end of the proof (see section \ref{testfunctions}), capitalizing on a new topological fact proved in section \ref{topological}.

\medskip \textbf{Acknowledgements} : I thank my thesis advisor Olivier Druet for stimulating discussions, his support, and his valuable remarks on a first draft of the paper. I would also like to thank Bruno Sévennec for his contribution in the decisive topological point (claim \ref{topology}).

\section{Construction of test functions} \label{construction}

In this section, we describe the construction of Nadirashvili \cite{NAD2002} (see also \cite{GIR2009}) which is at the basis of our theorem as well as of the previous results. Let $g$ be a metric on $\mathbf{S}^n$ conformal to $g_0$ of volume $1$. We denote by $dv_g$ the measure associated to $g$. We shall use in this paper the min-max characterization of the second eigenvalue of the Laplacian which tells us in particular that 
\begin{equation} \label{minmax} \lambda_2(\mathbf{S}^n,g) \leq \sup_{u\in E\setminus\{0\}} \frac{\int_{\mathbf{S}^n} \left|\nabla_g u\right|^2_g dv_g}{\int_{\mathbf{S}^n}  u^2 dv_g} 
\end{equation}
for all 2-dimensional subspaces $E$ of functions in $H^1\left(\mathbf{S}^n\right)$ with mean value $0$. The aim is to find a suitable space $E$ of test-functions such that (\ref{minmax}) gives the estimate of the theorem. 

\medskip On $\left({\mathbf S}^n,g_0\right)$, the eigenspace associated to $\lambda_1(\mathbf{S}^n,g_0)$ has dimension $n+1$ : it is the set of linear forms of $\mathbf{R}^{n+1}$ written $X_s = (s,.)$ for $s\in\mathbf{R}^{n+1}$. We will build $E$ with these functions, and as Hersch did for $\lambda_1(\mathbf{S}^n,g)$, we proceed to a renormalisation of measures in order to keep the orthogonality to constants. For $\xi\in\mathbf{B}^{n+1}$, we let $d_{\xi} : \overline{\mathbf{B}^{n+1}} \rightarrow \overline{\mathbf{B}^{n+1}}$ be defined  by 
$$ d_{\xi}(x) = \frac{(1-\left|\xi\right|^2)x +(1+2\xi.x +\left|x\right|^2)\xi}{1+2\xi.x +\left|x\right|^2\left|\xi\right|^2} $$
which is a conformal transformation when restricted to the unit sphere.

\medskip We say that $d_{\xi}$ renormalizes a finite measure $d\nu$ on the $n$-sphere if 
$$ \forall s\in \mathbf{S}^n, \int_{\mathbf{S}^n} X_s \circ d_{\xi}d\nu = 0 \hskip.1cm.$$
The Hersch lemma says that for all finite measures $d\nu$, such a $\xi$ exists. Moreover it is unique and depends continuously on $d{\nu}$ (the set of finite measures is considered as the topological dual of the continuous bounded functions) as proved in \cite{GIR2009}, Proposition 4.1.5. We call $\xi$ the renormalization point of $d\nu$.

\medskip We also define families of measures parametrized by the set of caps of $\mathbf{S}^n$, denoted by $\mathcal{C}$ :
$$a_{0,p} = \{ x\in \mathbf{S}^n ; x.p > 0 \} \hspace{20mm} a_{r,p} = d_{rp}\left(a_{0, p}\right) \hspace{5mm} (r,p)\in (-1,1)\times \mathbf{S}^n$$
We denote by $d\mu_a$ the "lift" of the measure $dv_g$ by the cap $a\in\mathcal{C}$ :
$${d\mu}_a = \left\{
    \begin{array}{ll}
        dv_g + \left(\tau_a\right)^* dv_g & \mbox{on } a \\
        0  & \mbox{on } a^*
    \end{array}
\right. $$

where $a^* = \mathbf{S}^n \setminus \bar{a}$ and $\tau_a$ is the conformal reflection with respect to the boundary circle of $a$, that is
$$\tau_{a_{r,p}} = d_{rp}\circ R_p \circ d_{-rp}$$
where 
$$R_p(x) = x - 2(p, x)p$$
is the reflection of $\mathbf{R}^{n+1}$ with respect to the hyperplane orthogonal to $p$. Let $\xi(a)$ be the renormalization point of $d\mu_a$. We set $d\nu_a = (d_{\xi(a)})_{*}d\mu_a$. Thanks to this family of measures, we can define a new family of test functions orthogonal to the constants~:
$$u^s_a = \left\{
    \begin{array}{ll}
        X_s \circ d_{\xi(a)} & \mbox{on } a \\
        X_s \circ d_{\xi(a)} \circ \tau_a & \mbox{on } a^*
    \end{array}
\right.$$
By a H\"older inequality, the numerator of the Rayleigh quotient is less than a conformal invariant.
\begin{equation} \label{gradient} 
\begin{aligned}
 \int_{\mathbf{S}^n} \left\vert \nabla_g u^s_a\right\vert^2_g dv_g 
 &< \left(\int_{\mathbf{S}^n} \left|\nabla_g u^s_a\right|^n_g dv_g \right)^{\frac{2}{n}} \\
 &= \left( 2 \int_{d_{\xi(a)}(a)} \left|\nabla_g X_s\right|^n_g dv_g \right)^{\frac{2}{n}} \\
 &< \left( 2\int_{\mathbf{S}^n} \left|\nabla_{g_0} X_s\right|^n_{g_0}dv_{g_0} \right)^{\frac{2}{n}}\\
 \end{aligned}
 \end{equation}
Let us define the multiplicity of a finite measure :

\begin{de} The multiplicity of a finite measure $d\nu$ on $\mathbf{S}^n$ is the dimension of the eigenspace $W$ associated to the maximal eigenvalue of the quadratic form :
$$Q(s) = \int_{\mathbf{S}^n} X_s^2 d\nu $$
We say that $d\nu$ is multiple if its multiplicity is greater than or equal to $2$. Otherwise, we say that $d\nu$ is simple.
\end{de}

As was noticed in \cite{GIR2009}, we know that if $dv_g$ is multiple, then we can choose $E = \{X_s; s\in W\}$ in (\ref{minmax}) to get that $\lambda_2(\mathbf{S}^n,g) \leq n(2\sigma_n)^{\frac{2}{n}}$. We also know that if there is a cap $a\in\mathcal{C}$ such that $d\nu_a$ is multiple, $\lambda_2(\mathbf{S}^n,g) < K_n n(2\sigma_n)^{\frac{2}{n}}$ using the space of test functions $E = \{u_a^s; s\in W\}$ in (\ref{minmax}). In this case, the theorem would be proved. In \cite{GIR2009}, it was proved that there necessarily exists such a multiple measure in odd dimensions (see below). 

\medskip Let us now assume that all measures $dv_g$ and $d\nu_a$, for $a\in\mathcal{C}$, are simple. Up to a renormalisation and a rotation, we may assume that 
$$ \forall t \in \mathbf{S}^n, \int_{\mathbf{S}^n} X_t dv_g = 0 $$
and that 
$$ \forall t \in \mathbf{S}^n \setminus [e_1], \int_{\mathbf{S}^n} X_t^2 dv_g < \int_{\mathbf{S}^n} X_{e_1}^2 dv_g \hskip.1cm.$$
We denote by $[s(a)]$ the unique direction of maximization of the quadratic form associated to $d\nu_a$. With the parametrization $(r,p)\in (-1,1)\times \mathbf{S}^n$ of $\mathcal{C}$, the maps $\xi : \mathcal{C} \rightarrow \mathbf{B}^{n+1}$ and $[s]:\mathcal{C} \rightarrow \mathbf{R}P^n$ are continuous. Moreover, one may prove that if $r\to -1$, that is $a\to\mathbf{S}^n$, we have :
\begin{equation} \label{limits} \lim_{a \to \mathbf{S}^n} \xi(a) = 0 \hspace{20mm} \lim_{a \to \mathbf{S}^n} [s(a)] = [e_1] \end{equation}

\section{Properties of the lift of the maximal direction}\label{topological}

Let us study the maps $\xi$ and $[s]$ at the light of the links between a cap $a\in\mathcal{C}$ and its symmetrical cap $a^* = \mathbf{S}^n \setminus \bar{a}$. With the parameter $(r,p)\in (-1,1)\times \mathbf{S}^n$, notice that $a_{r,p}^* = a_{-r,-p}$.

\begin{cl} \label{sym1} For $a\in\mathcal{C}$, we write $\xi^* = \xi(a^*)$, $[s^*] = [s(a^*)]$. Then 
$$ -\xi^{*} = \tau_a(-\xi) \hspace{10mm}\hbox{and}\hspace{10mm} [s^{*}] = R_a[s] $$
where $R_a = d_{\xi^*(a)} \circ \tau_a \circ d_{-\xi(a)}$ is an orthogonal map.
\end{cl}

\begin{pr} We set $\eta = -\tau_a(-\xi)$. Let $t\in\mathbf{S}^n$, then 
$$ \int_{\mathbf{S}^n} X_t \circ d_{\eta}\, d\mu_{a^{*}} = \int_{\mathbf{S}^n} X_t \circ d_{\eta} \circ \tau_a \, d\mu_{a} \hskip.1cm.$$
One can check that $d\mu_{a^{*}} = (\tau_a)^{*} d\mu_{a}$. The map $R_a = d_{\eta} \circ \tau_a \circ d_{-\xi(a)}$ is orthogonal because it is a M\"obius transformation of the unit ball preserving the origin (\cite{BEA1983}, Theorem 3.4.1). Thus we have that 
$$ \int_{\mathbf{S}^n} X_t \circ d_{\eta} \, d\mu_{a^{*}} = \int_{\mathbf{S}^n} X_t \circ R_a \circ d_{\xi}\, d\mu_{a} = \int_{\mathbf{S}^n} X_{R_a^{-1}(t)} \circ d_{\xi} \, d\mu_{a} = 0\hskip.1cm.$$
This is true for all $t\in\mathbf{S}^n$, and uniqueness of the renormalization point ensures that $\xi^*=\eta$.

\medskip The same argument with the function $\left(X_t \circ d_{\xi^*}\right)^2$ leads to 
$$ \forall t \in \mathbf{S}^n, \int_{\mathbf{S}^n} \left(X_t \circ d_{\xi^*}\right)^2 d\mu_{a^{*}} = \int_{\mathbf{S}^n} \left(X_{R_a^{-1}(t)} \circ d_{\xi}\right)^2 d\mu_{a} $$
and once again, we can conclude by uniqueness of the maximal direction that $[s^*] = R_a [s]$.
\end{pr}

\begin{rem} Thanks to this claim \ref{sym1}, we can prove the theorem in odd dimensions. Indeed, when $r\to 1$ that is $a\to\{p\}$, we use (\ref{limits}) in order to obtain :
$$ \lim_{a\to\{p\}}R_a = R_p $$
Then, $[s(a)] = R_a^{-1}[s^{*}(a)] \to R_p[e_1]$ when $a\to \{p\}$ by (\ref{limits}). Therefore, following~\cite{GIR2009} in odd dimensions, the map $[s] : [-1,1]\times \mathbf{S}^n \rightarrow \mathbf{R}P^n$ defines a homotopy between the constant map $[e_1]$ of degree $0$ and $\phi(p) = R_p[e_1]$ of degree $4$. Thus, there is a contradiction and there exists a multiple measure among $dv_g$ and $d\nu_a$ for $a\in\mathcal{C}$. 

\medskip We do not prove that the assumption that all measures are simple lead to a contradiction. Indeed, it is not clear that in even dimensions, such a configuration can not happen. Instead, we look for suitable test functions like in Nadirashvili's proof in dimension $2$~\cite{NAD2002}. However, inspired by the method of~\cite{GIR2009}, we use a topological argument to get symmetric properties of the lifts of the maximal directions. 
\end{rem}

\medskip The continuous map $[s] : [-1,1)\times \mathbf{S}^n \rightarrow \mathbf{R}P^n$ has exactly two continuous lifts because the set $[-1,1) \times \mathbf{S}^n$ is simply connected. We denote by $s$ the continuous lift such that $s(-1,.) = -e_1$, the other continuous lift is $-s$. Thanks to claim \ref{sym1},
$$ s(-r,-p) = \epsilon(r,p) R_{a_{r,p}} s(r,p) $$
where $\epsilon : [-1,1)\times \mathbf{S}^n \rightarrow \{\pm 1\}$ is a continuous map. Since $s\neq 0$ and $[-1,1)\times \mathbf{S}^n$ is connected, $\epsilon$ is a constant map. 

\medskip \begin{cl}\label{epsilon}
We have that $\epsilon=-1$. In other words, 
$$ s\left(a^{*}\right) = - R_{a} s(a) $$
for all caps $a$. 
\end{cl}

\begin{pr} We assume by contradiction that $\epsilon=1$. We set $f(p) = s(0,p)$ for $p\in\mathbf{S}^n$. This function $f$ is continuous on the sphere and satisfies 
\begin{equation} \label{sym2} \forall p \in \mathbf{S}^n, f(-p) = R_p f(p) \end{equation}
Indeed, $R_{a_{0,p}} =R_p$ because $\tau_{a_{0,p}} = R_p$. Using claim \ref{topology} below, we know that such a map $f$ can not have degree $0$. However, the map $s : [-1,0]\times\mathbf{S}^n \rightarrow \mathbf{S}^n$ defines a homotopy between $s_0 = f$ and $s_{-1} = -e_1$ of degree zero. Thus, there is a contradiction.
\end{pr}

\medskip We have used the following topology result~:

\begin{cl} \label{topology} Let $f : \mathbf{S}^n \rightarrow \mathbf{S}^n$ a continuous map which satisfies (\ref{sym2}). Then, if $n$ is odd, $\deg(f) = 1$ and if $n$ is even, $\deg(f) \in 2\mathbf{Z}+1$.
\end{cl}

\begin{pr} We first prove the claim for smooth functions which have a property of transversality (step 1) and we show that this case is generic (step 2).

\medskip {\it Step 1 - Let $f : \mathbf{S}^n \rightarrow \mathbf{S}^n$ be a smooth function which satisfies (\ref{sym2}). Let us assume that for all fix point $x \in \mathbf{S}^n$ of $f$, $T_x f - I : T_x\mathbf{S}^n \rightarrow T_x\mathbf{S}^n$ is an isomorphism. Then, if $n$ is odd, $\deg(f) = 1$ and if $n$ is even, $\deg(f) \in 2\mathbf{Z}+1$.}

\medskip {\it Proof of step 1 - }Let $F$ be defined by 
$$\begin{array}[t]{lrcl}
F : & \mathbf{S}^n \times [-1,1] & \longrightarrow & \mathbf{R}^{n+1} \\
    & (x,t) & \longmapsto & \frac{1}{2}\left( f(x) - x + t (f(x)+x) \right)
    \end{array}$$
We notice that if $F$ never vanishes, $\frac{F}{\left|F\right|}$ defines a homotopy between $f$ and $\sigma$, the antipodal map and $\deg(f) = \deg(\sigma) = (-1)^{n+1}$.

Now, $F(x,t) = 0$ if and only if $t=0$ and $x$ is a fix point of $f$ and then,
$$ \forall (v,t) \in T_x\mathbf{S}^n \times \mathbf{R}, \hspace{5mm} DF(x,0)(v,t) = \frac{1}{2}(T_x f - I)v + x t \hskip.1cm.$$
Thus, $DF(x,0)$ is an isomorphism, and $0$ is a regular value. We write $(x_1,0),\cdots,(x_r,0)$ the regular points of $F^{-1}(0)$. Let's approximate $F$ by its differential in the neighborhood of its zeros. Let $\alpha >0$ and, set for $1\leq i \leq r$, $\phi_i : B_{x_i}(\alpha) \rightarrow B_0(\alpha) \subset T_{x_i} \mathbf{S}^n$ the exponential chart at $x_i$. We obtain for $(x,t) \in B_{x_i}(\alpha)\times (-\alpha,\alpha)$ 
$$ F(x,t) = DF(x_i,0)(\phi_i(x),t) + R_i(\phi_i(x),t) $$
where $\frac{R_i(v,t)}{\left|(v,t)\right|} \to 0$ when $(v,t)\to 0$. We write for $x\in\mathbf{S}^n$ that 
$$ F_t(x) = F(x,t) \hspace{15mm} L_t(x) = \left\{
    \begin{array}{ll}
        DF(x_i,0)(\phi_i(x),t) & \mbox{if } (x,t) \in B_{x_i}(\alpha)\times (-\alpha,\alpha) \\
        0 & \mbox{otherwise.}
    \end{array}
\right.
$$
We define a cut-off function $0 \leq \psi \leq 1$ such that $\psi = 1$ on $K_1 = \bigcup_{i=1}^{r} \overline{B_{x_i}(\frac{\alpha}{2})}$ and $\psi=0$ on $K_2 = \mathbf{S}^n \setminus \bigcup_{i=1}^{r} B_{x_i}(\alpha)$. We set for $s\in[0,1]$ 
$$ G_s^t = \frac{s\psi L_t + (1-s\psi) F_t}{\left|s\psi L_t + (1-s\psi) F_t\right|} \hskip.1cm.$$
One may choose $\alpha>0$ small enough so that $G_{s}^t$ is well defined for all $t \in (-\alpha,\alpha)\setminus \{0\}$. Then, for $0< t< \alpha $, $G_1^t$ is homotopic to $G_0^t = \frac{F_t}{\left|F_t\right|}$, so to $f$, and $G_1^{-t}$ is homotopic to $\sigma$. We now write, for $t\in(-\alpha,\alpha)$, $g_t = G_1^t$. 

Let us look at the behaviour of $g_t = \frac{L_t}{\left|L_t\right|}$ in the balls $\overline{B_{x_i}(\frac{\alpha}{2})}$ when $t\to 0$. We recall that 
$$L_t(x) = \frac{1}{2}(T_{x_i} f - I) \phi_i(x) + x_i t \hskip.1cm.$$
Therefore, the image $I_{x_i}^t = g_t(\overline{B_{x_i}(\frac{\alpha}{2})})$ blows up to the half-sphere $D_{x_i} = \{x\in\mathbf{S}^n ; (x,x_i) >0\}$ when $t\to 0$. 

Thanks to (\ref{sym2}), $x$ is a fix point of $f$ if and only if $-x$ is a fix point too. Moreover, by differentiating (\ref{sym2}) at a fix point $x$, we obtain $T_{-x} f - I = - (T_x f -I)$. 

Let's renumber the fix points $x_1,\cdots,x_k,-x_1,\cdots,-x_k$ (with $r=2k$), so that $x_1,\cdots,x_k$ are in a same half sphere $D_p = \{(x,p) >0 \}$. We choose $\epsilon <\alpha$ small enough so that $\bigcap_{i=1}^{k} I_{x_i}^{\epsilon}$ has a non-empty interior $I$. Then, for $z\in I$, there is a unique point in $g_t^{-1}(z) \cap B_{x_i}(\frac{\alpha}{2})$ for all $0 < t < \epsilon$. Since $g_{\epsilon}(x) = g_{-\epsilon}(-x)$, if $z \in I$, then $z \in I_{-x_i}^{-\epsilon}$ and $z\notin I_{-x_i}^{\epsilon} \cup I_{x_i}^{-\epsilon}$. 

For $1\leq i \leq k$, let $\{a_i\} = B_{x_i}(\frac{\alpha}{2}) \cap g_{\epsilon}^{-1}(z)$. Then by definition of degree and homotopy,
$$ \deg(f) - \deg(\sigma) = \deg(g_{\epsilon}) - \deg(g_{-\epsilon}) = \sum_{i=1}^{k} \operatorname{ind}_{a_i}(g_{\epsilon}) - \operatorname{ind}_{-a_i}(g_{-\epsilon}) = \sum_{i=1}^{k} (1-(-1)^{n+1}) \nu_i $$
where $\nu_i = \operatorname{ind}_{a_i}(g_{\epsilon}) \in \pm 1$. In odd dimensions, $\deg(f)= \deg(\sigma) = 1$ and in even dimensions, $\deg(f) \in 2\mathbf{Z}+1$. This ends the proof of step 1.

\medskip {\it Step 2 - Let $f : \mathbf{S}^n \rightarrow \mathbf{S}^n$ be a continuous map which satisfies (\ref{sym2}). Then there exists a map, homotopic to $f$, which satisfies the assumptions of step 1.}

\medskip {\it Proof of step 2 -} Denote by $(e_0,e_1,\cdots,e_n)$ the canonical basis of $\mathbf{R}^{n+1}$ and $B_k^{\alpha} \subset D_{e_k} = \{(x,e_k)>0\}$ the ball centered at $e_k$ such that $d(B_k^{\alpha},D_{-e_k}) = \alpha >0$. Choose $\alpha$ small enough so that 
$$\bigcup_{i=0}^{n} B_i^{2\alpha} \cup (-B_i^{2\alpha}) = \mathbf{S}^n \hskip.1cm.$$
Let $\epsilon >0$. We build by induction maps $g_k : \mathbf{S}^n \rightarrow \mathbf{S}^n$ such that $g_0 = f$ and, for $0\leq k\leq n$, 
\begin{itemize}
\item $g_{k+1} = g_k$ on $\mathbf{S}^n \setminus \left(B_k^{\alpha} \cup (-B_k^{\alpha})\right)$
\item $g_{k+1}$ is smooth on $\bigcup_{i=0}^{k} B_i^{2\alpha} \cup (-B_i^{2\alpha})$
\item $\left\|g_{k+1}-g_k\right\|_{\mathcal{C}^0} <\epsilon$
\item $g_{k+1}$ satisfies ($\ref{sym2}$).
\end{itemize}
By density of smooth maps $\mathbf{S}^n \rightarrow \mathbf{R}^{n+1}$, choose $h_k$ such that $\left\|h_k - g_k\right\|_{\mathcal{C}^0} < \epsilon$. Let $0\leq \phi\leq 1$ be a smooth cut-off function such that $\phi = 1$ on $B_i^{2\alpha}$ and $\phi = 0$ on $\mathbf{S}^n\setminus B_i^{\alpha}$. We let $g_{k+1}$ be defined, provided $\epsilon$ is small enough, by 
$$ g_{k+1}(x) = \frac{\phi h_k + (1-\phi)g_k}{\left|\phi h_k + (1-\phi)g_k\right|} \hbox{ and }g_{k+1}(-x) = R_x \circ g_{k+1}(x) $$
for $x\in \overline{D_{e_k}}$. Therefore $g = g_{n+1}$ is smooth, satisfies ($\ref{sym2}$) and $\left\|g-f\right\|_{\mathcal{C}^0} < C\epsilon$. If $\epsilon$ is small enough, $g$ is homotopic to $f$. 

Let's now tackle the transversality condition. We write $g$ in the following way 
$$ g(x) = X(x) + \lambda(x)x $$
where $X$ is a tangent vector field of the sphere and $\left|X\right|^2 + \lambda^2 = 1$. Then, $g$ satisfies ($\ref{sym2}$) if and only if $X$ and $\lambda$ are even maps. By differentiating these equalities at a fix point $x$ (with $\lambda(x) = 1$ and $X(x)=0$), one may find $T_x g - I = T_x X$. Then, $T_x g - I$ is an isomorphism for all fix points $x$ if and only if $X$ is transverse to the zero vector field. Then, one may build by induction, with Sard's theorem in $n$-dimensional charts on $D_{e_k}$, smooth tangent vector fields $X_k$ such that $X_0 = X$ and for $0\leq k \leq n$ :
\begin{itemize}
\item $X_{k+1} = X_k$ on $\mathbf{S}^n \setminus \left(B_k^{\alpha} \cup (-B_k^{\alpha})\right)$
\item $X_{k+1}$ is transverse to $0$ on $\bigcup_{i=0}^{k} B_i^{2\alpha} \cup (-B_i^{2\alpha})$
\item $\left\|X_{k+1}-X_k\right\|_{\mathcal{C}^0} <\epsilon$
\item $X_{k+1}$ is an even map.
\end{itemize}
Set $\bar{f}(x) = \frac{X_{n+1}(x) + \lambda(x)x}{\left|X_{n+1}(x)\right|^2 + \lambda(x)^2}$. If $\epsilon$ is small enough, then $\bar{f}$ is well defined, satisfies the assumptions of step 1 and is homotopic to $f$. This ends the proof of step 2.

\medskip These two steps clearly end the proof of the claim. 
\end{pr}

\section{Choice of test functions} \label{testfunctions}
Thanks to claim \ref{epsilon}, one may easily deduce that 
\begin{equation} \label{sym4} \forall a\in\mathcal{C}, {u}_{a^*} = - u_a \end{equation}
where we have set, for this section $u_a=u_a^{s(a)}$. Let $r\in (-1,1)$. We look at the space $E$ generated by
$$ \phi = X_{e_1} \hbox{ and } \psi_r = u_{a_{r,e_1}} \hskip.1cm.$$
One may deduce from the continuity of $\xi$ and $s$, (\ref{limits}) and (\ref{sym4}), that

\begin{cl} \label{convergence} The map $r\in(-1,1) \mapsto \psi_r \in \left(L^2(\mathbf{S}^n,g), \left\|.\right\|_{L^2} \right)$ is continuous and 
$$ \lim_{r \to -1} \psi_r = -\phi \hspace{20mm} \lim_{r \to 1} \psi_r = \phi$$
\end{cl}

For $(x,y)\in \mathbf{R}^2\setminus\{0\}$, we set $f_r = x\phi + y\psi_r \in E$. Conformal invariance gives that 
$$ \frac{\int_{\mathbf{S}^n}\left|\nabla_g f_r\right|^2_g dv_g}{\int_{\mathbf{S}^2} f_r^2 dv_g} = \frac{C_n}{2^{\frac{2}{n}}} \frac{\sigma x^2 + \tau_r y^2 + 2 \alpha_r xy}{I x^2 + J_r y^2 + 2 \beta_r xy} := \frac{C_n}{2^{\frac{2}{n}}} q(x,y) $$
where $(n+1)\left(\int_{\mathbf{S}^n}\left|\nabla_g \phi \right|_g^n dv_g\right)^{\frac{2}{n}} = C_n = K_n n(2\sigma_n)^{\frac{2}{n}}$ and we set for $r\in(-1,1)$ 
$$\sigma = \frac{\int_{\mathbf{S}^n}\left|\nabla_g \phi \right|_g^2 dv_g}{\left(\int_{\mathbf{S}^n}\left|\nabla_g \phi \right|_g^n dv_g\right)^{\frac{2}{n}}} < 1 \hspace{15mm} \tau_r = \frac{\int_{\mathbf{S}^n}\left|\nabla_g \psi_r \right|_g^2 dv_g}{\left(\int_{\mathbf{S}^n}\left|\nabla_g \phi \right|_g^n dv_g\right)^{\frac{2}{n}}} < 2^{\frac{2}{n}}$$
$$ \alpha_r = \frac{\int_{\mathbf{S}^n}g(\nabla_g \psi_r,\nabla_g \phi) dv_g}{\left(\int_{\mathbf{S}^n}\left|\nabla_g \phi \right|_g^n dv_g\right)^{\frac{2}{n}}} \hspace{15mm} \beta_r = (n+1)\int_{\mathbf{S}^n} \phi\psi_r dv_g $$
$$ I = (n+1)\int_{\mathbf{S}^n} \phi^2 dv_g > 1 \hspace{15mm} J_r = (n+1)\int_{\mathbf{S}^n} \psi_r^2 dv_g > 1$$
By (\ref{gradient}), $\tau_r < 2^{\frac{2}{n}}$ and by maximality of $\phi$ and $\psi_r$, $I >1$ and $J_r>1$.

\noindent Thus, in order to get the estimate of the theorem and using the min-max principle (\ref{minmax}), we look for $r\in(-1,1)$ such that for all $(x,y) \in \mathbf{R}^2 \setminus \{0\}$ :
$$ q(x,y) < 2^{\frac{2}{n}} \hskip.1cm. $$
Since $I>1$ and $J_r>1$, we look for $r\in(-1,1)$ such that 
$$ (\sigma-2^{\frac{2}{n}})x^2 + 2(\alpha_r - 2^{\frac{2}{n}}\beta_r) yx + (\tau_r - 2^{\frac{2}{n}}) y^2 < 0 \hskip.1cm.$$
Moreover, since $\sigma < 1$ and $\tau_r - 2^{\frac{2}{n}} < 0$, it is sufficient to find $r\in(-1,1)$ such that 
$$\alpha_r - 2^{\frac{2}{n}}\beta_r = 0\hskip.1cm.$$
By the claim \ref{convergence}, we know that 
$$ \alpha_r = \frac{ - \int_{\mathbf{S}^n}\psi_r \left(\Delta_g \phi\right) dv_g}{\left(\int_{\mathbf{S}^n}\left|\nabla_g \phi \right|_g^n dv_g\right)^{\frac{2}{n}}} \underset{r\to 1}{\longrightarrow} \frac{ - \int_{\mathbf{S}^n}\phi \left(\Delta_g \phi\right) dv_g}{\left(\int_{\mathbf{S}^n}\left|\nabla_g \phi \right|_g^n dv_g\right)^{\frac{2}{n}}} = \sigma $$
and that 
$$ \beta_r = (n+1)\int_{\mathbf{S}^n} \phi\psi_r dv_g \underset{r\to 1}{\longrightarrow} (n+1)\int_{\mathbf{S}^n} \phi^2 dv_g  = I\hskip.1cm. $$
Thus, when $r\to 1$ and in an analogous way, when $r\to -1$, (see claim \ref{convergence}),
$$ \alpha_r - 2^{\frac{2}{n}}\beta_r \underset{r\to 1}{\longrightarrow} \sigma - 2^{\frac{2}{n}}I < 0 $$
and
$$ \alpha_r - 2^{\frac{2}{n}}\beta_r \underset{r\to -1}{\longrightarrow} 2^{\frac{2}{n}}I - \sigma > 0 \hskip.1cm.$$
By continuity, (claim \ref{convergence}), there exists $r\in (-1,1)$ such that $\alpha_r - 2^{\frac{2}{n}}\beta_r = 0$. As already said, this completes the proof of the theorem.

\bibliographystyle{plain}
\bibliography{biblio}
\end{document}